\documentclass[11pt,reqno]{amsart}
\usepackage{amsfonts,amssymb,amsmath}
\usepackage{color}

\DeclareMathOperator{\lin}{\mathrm{span}}

\newtheorem{theorem}{Theorem}

\newtheorem{proposition}[theorem]{Proposition}

\numberwithin{equation}{section} \numberwithin{theorem}{section}

\begin{document}

\title{Pattern Recognition on Oriented Matroids:
Decompositions of Topes,\\ and Orthogonality Relations}

\author{Andrey O. Matveev}
\email{andrey.o.matveev@gmail.com}



\begin{abstract}
If $\mathrm{V}(\boldsymbol{R})$ is the vertex set of a symmetric cycle $\boldsymbol{R}$ in the tope graph of a simple oriented matroid $\mathcal{M}$, then for any tope $T$ of $\mathcal{M}$ there exists a unique inclusion-minimal subset $\boldsymbol{Q}(T,\boldsymbol{R})$ of $\mathrm{V}(\boldsymbol{R})$ such that $T$ is the sum of the topes of $\boldsymbol{Q}(T,\boldsymbol{R})$.

If for decompositions $\boldsymbol{Q}(T',\boldsymbol{R}')$ and $\boldsymbol{Q}(T'',\boldsymbol{R}'')$
with respect to symmetric cycles $\boldsymbol{R}'$ and $\boldsymbol{R}''$ in the tope graphs of two simple oriented matroids, whose ground sets have the cardinalities of opposite parity, we have $|\boldsymbol{Q}(T',\boldsymbol{R}')|\geq 5$ and $|\boldsymbol{Q}(T'',\boldsymbol{R}'')|\geq 5$, then these decompositions satisfy a certain orthogonality relation.
\end{abstract}

\maketitle

\pagestyle{myheadings}

\markboth{A.O.~MATVEEV}{PATTERN RECOGNITION ON ORIENTED MATROIDS}

\thispagestyle{empty}


\vspace{-3mm}
\section{Introduction}

Let $\mathcal{M}':=(E_s,\mathcal{L}')=(E_s,\mathcal{T}')$ and $\mathcal{M}'':=(E_t,\mathcal{L}'')=(E_t,\mathcal{T}'')$ be {\em simple\/} (i.e., with no {\em loops}, {\em parallel elements\/} or {\em antiparallel elements}) oriented matroids on their ground sets $E_s=[s]:=\{1,\ldots,s\}$ and $E_t$, with sets of covectors~$\mathcal{L}'$ and~$\mathcal{L}''$, and with sets of topes~$\mathcal{T}'\subseteq\{1,-1\}^{E_s}$ and~$\mathcal{T}''\subseteq\{1,-1\}^{E_t}$, respectively; see~\cite{BLSWZ} on oriented matroids. We suppose that
\begin{gather}
\nonumber
s<t\; ,\\
\label{prchap:10:eq:5}
s\not\equiv t\pmod{2}\; .
\end{gather}

Let $\boldsymbol{R}'$ be a {\em symmetric $2s$-cycle}, with its vertex set $\mathrm{V}(\boldsymbol{R}')=-\mathrm{V}(\boldsymbol{R}')$, in the tope graph~$\mathcal{T}'(\mathcal{L}')$ of~$\mathcal{M}'$, and let~$\boldsymbol{R}''$ be a {\em symmetric~$2t$-cycle\/} in the tope graph~$\mathcal{T}''(\mathcal{L}'')$ of~$\mathcal{M}''$.

Let $T'\in\mathcal{T}'$ and $T''\in\mathcal{T}''$ be any topes of the oriented matroids~$\mathcal{M}'$ and~$\mathcal{M}''$ such that for the {\em unique inclusion-minimal subsets\/}~$\boldsymbol{Q}(T',\boldsymbol{R}')\subset\mathrm{V}(\boldsymbol{R}')$ and~$\boldsymbol{Q}(T'',\boldsymbol{R}'')\subset\mathrm{V}(\boldsymbol{R}'')$ with the properties
\begin{align*}
\sum_{Q'\in\boldsymbol{Q}(T',\boldsymbol{R}')}Q'=T'\ \ \ &\text{and} \ \ \ \sum_{Q''\in\boldsymbol{Q}(T'',\boldsymbol{R}'')}Q''=T''\; ,
\intertext{discussed in~\cite[Sect.~3.11]{M-DG-I}, we have}
|\boldsymbol{Q}(T',\boldsymbol{R}')|\geq 5\ \ \ &\text{and} \ \ \ |\boldsymbol{Q}(T'',\boldsymbol{R}'')|\geq 5\; .
\end{align*}

Let
\begin{equation*}
\boldsymbol{\Lambda}'':=\boldsymbol{\Lambda}(T'',\boldsymbol{R}'')
\end{equation*}
be the abstract simplicial complex on the vertex set $E_t$, with the facet family
\begin{equation}
\label{prchap:10:eq:9}
\bigl\{E_t-\mathbf{S}(T'',Q'')\colon Q''\in\boldsymbol{Q}(T'',\boldsymbol{R}'')\bigr\}\; ,
\end{equation}
where $\mathbf{S}(T'',Q'')$ is the {\em separation set\/} of the topes~$T''$ and $Q''$; see~e.g.~\cite{M-9} on such complexes.

Associate with the complex~$\boldsymbol{\Lambda}''$ its ``{\em long}'' {\em $f$-vector}
\begin{equation*}
\boldsymbol{f}(\boldsymbol{\Lambda}'';t):=\bigl(f_0(\boldsymbol{\Lambda}'';t),f_1(\boldsymbol{\Lambda}'';t),
\ldots,f_t(\boldsymbol{\Lambda}'';t)\bigr)\in\mathbb{N}^{t+1}
\end{equation*}
defined by
\begin{equation*}
f_j(\boldsymbol{\Lambda}'';t):=\#\{F\in\boldsymbol{\Lambda}''\colon |F|=j\}\; ,\ \ \ 0\leq j\leq t\; ;
\end{equation*}
throughout the paper, the components of all vectors, as well as the rows and columns of matrices are indexed starting with zero.

Define vectors $\boldsymbol{\beta}(t;t)\in\mathbb{P}^{t+1}$
and~$\boldsymbol{\beta}(s;t)\in\mathbb{N}^{t+1}$ by
\begin{align*}
\boldsymbol{\beta}(t;t):\!&=(\tbinom{t}{0},\tbinom{t}{1},\ldots,\tbinom{t}{t})\; ,\\
\boldsymbol{\beta}(s;t):\!&=(\tbinom{s}{0},\tbinom{s}{1},\ldots,\tbinom{s}{t})\; .
\end{align*}
Let $\mathbf{U}(t)$ denote the square {\em backward identity matrix\/} of order $t+1$, whose $(i,j)$th entry is the~Kronecker delta~$\delta_{i+j,t}$. We denote by $\mathbf{T}(t)$ the square {\em forward shift matrix\/} of order~$t+1$, whose $(i,j)$th entry is $\delta_{j-i,1}$.

According to the argument given in~\cite{M-9}, and by~\cite[Prop.~3.51(a)]{G-DG}, the vector
\begin{equation*}
\bigl(\boldsymbol{\beta}(t;t)-\boldsymbol{f}(\boldsymbol{\Lambda}'';t)\bigr)\mathbf{U}(t)\in\mathbb{N}^{t+1}
\end{equation*}
is the {\em long $f$-vector\/}~$\boldsymbol{f}(\boldsymbol{\Omega}'';t)$
of the {\em boundary complex\/}~$\boldsymbol{\Omega}''$ of a~$(t-3)$-dimensional {\em simplicial\/} convex {\em polytope\/} with $t$ vertices.

Associate with the abstract simplicial complex
\begin{equation*}
\boldsymbol{\Lambda}':=\boldsymbol{\Lambda}(T',\boldsymbol{R}')
\end{equation*}
on the vertex set $E_s$, with the facet family
\begin{equation}
\label{prchap:10:eq:10}
\bigl\{E_s-\mathbf{S}(T',Q')\colon Q'\in\boldsymbol{Q}(T',\boldsymbol{R}')\bigr\}\; ,
\end{equation}
its {\em long $f$-vector}
\begin{equation*}
\boldsymbol{f}(\boldsymbol{\Lambda}';t)\in\mathbb{N}^{t+1}
\end{equation*}
defined by $f_j(\boldsymbol{\Lambda}';t):=\#\{F\in\boldsymbol{\Lambda}'\colon |F|=j\}$, $0\leq j\leq t$.

The vector
\begin{equation*}
\bigl(\boldsymbol{\beta}(s;t)-\boldsymbol{f}(\boldsymbol{\Lambda}';t)\bigr)\mathbf{T}(t)^{t-s}\mathbf{U}(t)\in\mathbb{N}^{t+1}
\end{equation*}
is the {\em long $f$-vector\/}~$\boldsymbol{f}(\boldsymbol{\Omega}';t)$
of the {\em boundary complex\/}~$\boldsymbol{\Omega}'$ of an~$(s-3)$-dimensional {\em simplicial\/} convex {\em polytope\/} with $s$ vertices.

Define the $(i,j)$th entry of a square matrix
\begin{equation*}
\mathbf{S}(t)
\end{equation*}
of order~$t+1$ to be $(-1)^{j-i}\tbinom{t-i}{j-i}$.

Recall that the standard {\em $h$-vectors} $\boldsymbol{h}(\boldsymbol{\Omega}'')\in\mathbb{N}^{t-2}$ and~$\boldsymbol{h}(\boldsymbol{\Omega}')\in\mathbb{N}^{s-2}$ of the {\em boundary complexes\/} $\boldsymbol{\Omega}''$ and $\boldsymbol{\Omega}'$ of {\em simplicial polytopes\/} both satisfy the~{\em Dehn--Sommerville relations\/} \cite[Sect.~8.3]{Z}
\begin{align*}
h_k(\boldsymbol{\Omega}'')&=h_{t-k-3}(\boldsymbol{\Omega}'')\; , & 0\leq&\; k\leq t-3\; ;\\
h_k(\boldsymbol{\Omega}')&=h_{s-k-3}(\boldsymbol{\Omega}')\; , & 0\leq&\; k\leq s-3\; .
\end{align*}
As a consequence, the ``{\em long}'' {\em $h$-vectors}
\begin{equation}
\label{prchap:10:eq:1}
\boldsymbol{h}(\boldsymbol{\Omega}'';t):=\bigl(\boldsymbol{\beta}(t;t)-\boldsymbol{f}(\boldsymbol{\Lambda}'';t)\bigr)\mathbf{U}(t)\mathbf{S}(t)\in\mathbb{Z}^{t+1}
\end{equation}
and
\begin{equation}
\label{prchap:10:eq:2}
\boldsymbol{h}(\boldsymbol{\Omega}';t):=\bigl(\boldsymbol{\beta}(s;t)-\boldsymbol{f}(\boldsymbol{\Lambda}';t)\bigr)\mathbf{T}(t)^{t-s}\mathbf{U}(t)\mathbf{S}(t)
\in\mathbb{Z}^{t+1}
\end{equation}
of the complexes~$\boldsymbol{\Omega}''$ and $\boldsymbol{\Omega}'$ satisfy the {\em Dehn--Sommerville\/} type {\em relations}
\begin{align}
\label{prchap:10:eq:7}
h_k(\boldsymbol{\Omega}'';t)&=-h_{t-k}(\boldsymbol{\Omega}'';t)\; , & 0\leq&\; k\leq t\; ;\\
\label{prchap:10:eq:8}
h_k(\boldsymbol{\Omega}';t)&=\phantom{-}h_{t-k}(\boldsymbol{\Omega}';t)\; , & 0\leq&\; k\leq t\; ;
\end{align}
see e.g.~\cite[Sect.~2.3]{M-DG-I}.

\section{Orthogonality relations for decompositions of topes with respect to symmetric cycles in the tope graphs}

Since the maximal face $[t]$ of the simplex $\mathbf{2}^{[t]}$ does not belong to the complexes~$\boldsymbol{\Omega}'$ and~$\boldsymbol{\Omega}''$, by~\cite[Eq.~(2.3) of Prop.~2.1]{M-DG-I} we have
\begin{equation*}
\bigl\langle\boldsymbol{h}(\boldsymbol{\Omega}';t),\boldsymbol{\iota}(t)\bigr\rangle=
\bigl\langle\boldsymbol{h}(\boldsymbol{\Omega}'';t),\boldsymbol{\iota}(t)\bigr\rangle=0\; ,
\end{equation*}
for the vector
\begin{equation*}
\boldsymbol{\iota}(t):=(1,1,\ldots,1)\in\mathbb{N}^{t+1}\; ,
\end{equation*}
where $\langle \cdot,\cdot\rangle$ means the standard scalar product.

For the positive integers $k$, define abstract simplicial complexes $\overline{\mathbf{2}^{[k]}}$ with their {\em long $h$-vectors\/}
$\boldsymbol{h}(\overline{\mathbf{2}^{[k]}};t):=\boldsymbol{f}(\overline{\mathbf{2}^{[k]}};t)\mathbf{S}(t)$
to be the {\em boundary complexes\/} of the {\em simplices\/}~$\mathbf{2}^{[k]}$ by
\begin{equation*}
\overline{\mathbf{2}^{[k]}}:=\mathbf{2}^{[k]}-\{[k]\}\; .
\end{equation*}

In view of~(\ref{prchap:10:eq:8}), the long $h$-vector $\boldsymbol{h}(\boldsymbol{\Omega}')$ lies either in the linear span
\begin{equation}
\label{prchap:10:eq:3}
\lin\Bigl(\boldsymbol{h}(\overline{\mathbf{2}^{[1]}};t),\boldsymbol{h}(\overline{\mathbf{2}^{[3]}};t),\ldots,\boldsymbol{h}(\overline{\mathbf{2}^{[s-2]}};t)\Bigr)\; ,
\end{equation}
when $t$ is {\em even}, or in the linear span
\begin{equation}
\label{prchap:10:eq:4}
\lin\Bigl(\boldsymbol{h}(\overline{\mathbf{2}^{[2]}};t),\boldsymbol{h}(\overline{\mathbf{2}^{[4]}};t),\ldots,\boldsymbol{h}(\overline{\mathbf{2}^{[s-2]}};t)\Bigr)\; ,
\end{equation}
when $t$ is {\em odd}; see~\cite[Sect.~3.1]{M-DG-I}.

The Dehn--Sommerville type relations~(\ref{prchap:10:eq:7}) and~(\ref{prchap:10:eq:8}) imply that~$\boldsymbol{h}(\boldsymbol{\Omega}';t)$ is a {\em left eigenvector\/} of the backward identity matrix~$\mathbf{U}(t)$ that corresponds to its eigenvalue~$1$, while $\boldsymbol{h}(\boldsymbol{\Omega}'';t)$
is a {\em right eigenvector\/} of $\mathbf{U}(t)$ that corresponds to the other eigenvalue~$-1$. By the {\em principle of biorthogonality\/}~\cite[Th.~1.4.7(a)]{HJ} we have
\begin{equation}
\label{prchap:10:eq:6}
\bigl\langle\boldsymbol{h}(\boldsymbol{\Omega}';t),\boldsymbol{h}(\boldsymbol{\Omega}'';t)\bigr\rangle=0\; .
\end{equation}
In other words, together with the relation~(\ref{prchap:10:eq:6}), the definitions~(\ref{prchap:10:eq:1}) and~(\ref{prchap:10:eq:2}) yield
\begin{equation*}
\bigl(\boldsymbol{\beta}(s;t)-\boldsymbol{f}(\boldsymbol{\Lambda}';t)\bigr)\mathbf{T}(t)^{t-s}
\mathbf{U}(t)\mathbf{S}(t)\mathbf{S}(t)^{\!\top}\mathbf{U}(t)\bigl(\boldsymbol{\beta}(t;t)^{\!\top}-\boldsymbol{f}(\boldsymbol{\Lambda}'';t)^{\!\top}\bigr)=0\; .
\end{equation*}
Note that the $(i,j)$th entry of the square matrix
\begin{equation*}
\mathbf{M}(t):=\mathbf{U}(t)\mathbf{S}(t)\mathbf{S}(t)^{\!\top}\mathbf{U}(t)
\end{equation*}
of order $t+1$ is $(-1)^{i+j}\tbinom{i+j}{i}$.

Let us sum up our conclusions:
\begin{proposition}
Let $\mathcal{M}':=(E_s,\mathcal{L}')=(E_s,\mathcal{T}')$ and $\mathcal{M}'':=(E_t,\mathcal{L}'')=(E_t,\mathcal{T}'')$ be simple oriented matroids on their ground sets $E_s$ and $E_t$ such that
\begin{gather*}
s<t\; ,\\
s\not\equiv t\pmod{2}\; ,
\end{gather*}
with sets of covectors~$\mathcal{L}'$ and~$\mathcal{L}''$, and with sets of topes~$\mathcal{T}'\subseteq\{1,-1\}^{E_s}$ and~$\mathcal{T}''\subseteq\{1,-1\}^{E_t}$, respectively.

Let $\boldsymbol{R}'$ be a symmetric cycle in the tope graph~$\mathcal{T}'(\mathcal{L}')$ of~$\mathcal{M}'$, and let~$\boldsymbol{R}''$ be a symmetric~cycle in the tope graph~$\mathcal{T}''(\mathcal{L}'')$ of~$\mathcal{M}''$.

Let $T'\in\mathcal{T}'$ and $T''\in\mathcal{T}''$ be any topes of the oriented matroids~$\mathcal{M}'$ and~$\mathcal{M}''$ such that for the unique inclusion-minimal subsets~$\boldsymbol{Q}(T',\boldsymbol{R}')\subset\mathrm{V}(\boldsymbol{R}')$ and~$\boldsymbol{Q}(T'',\boldsymbol{R}'')\subset\mathrm{V}(\boldsymbol{R}'')$ with the properties
\begin{align*}
\sum_{Q'\in\boldsymbol{Q}(T',\boldsymbol{R}')}Q'=T'\ \ \ &\text{and} \ \ \ \sum_{Q''\in\boldsymbol{Q}(T'',\boldsymbol{R}'')}Q''=T''\; ,
\intertext{we have}
|\boldsymbol{Q}(T',\boldsymbol{R}')|\geq 5\ \ \ &\text{and} \ \ \ |\boldsymbol{Q}(T'',\boldsymbol{R}'')|\geq 5\; .
\end{align*}

The long $f$-vectors $\boldsymbol{f}(\boldsymbol{\Lambda}'';t)$ and~$\boldsymbol{f}(\boldsymbol{\Lambda}';t)$ of the complexes~$\boldsymbol{\Lambda}''$
and~$\boldsymbol{\Lambda}'$ whose families of facets are defined by~{\rm(\ref{prchap:10:eq:9})} and~{\rm(\ref{prchap:10:eq:10})}, respectively, satisfy the {\em orthogonality relation}
\begin{equation*}
\bigl(\boldsymbol{\beta}(s;t)-\boldsymbol{f}(\boldsymbol{\Lambda}';t)\bigr)\mathbf{T}(t)^{t-s}\cdot
\mathbf{M}(t)\cdot\bigl(\boldsymbol{\beta}(t;t)^{\!\top}-\boldsymbol{f}(\boldsymbol{\Lambda}'';t)^{\!\top}\bigr)=0\; .
\end{equation*}
\end{proposition}

\vspace{3mm}


\begin{thebibliography}{10}
\bibitem{BLSWZ}
{\em Bj\"{o}rner A.}, {\em Las~Vergnas M.}, {\em Sturmfels B.}, {\em White N.}, {\em Ziegler G.M.} Oriented matroids. Second
edition. Encyclopedia of Mathematics, 46. -- Cambridge: Cambridge University Press, 1999.

\bibitem{G-DG}
{\em Gainanov D.N.} Graphs for pattern recognition. Infeasible systems of linear inequalities. -- Berlin: De Gruyter, 2016.

\bibitem{HJ}
{\em Horn R.A.}, {\em Johnson C.R.} Matrix analysis. Second edition. -- Cambridge: Cambridge University Press, 2013.

\bibitem{M-DG-I}
{\em Matveev A.O.} Pattern recognition on oriented matroids. -- Berlin: De Gruyter, 2017.

\bibitem{M-9}
{\em Matveev A.O.} Pattern recognition on oriented matroids: Decompositions of topes, and Dehn--Sommerville type relations. Preprint [arXiv:1703.04508],
2017.

\bibitem{Z}
{\em Ziegler G.M.} Lectures on polytopes. Revised edition. Graduate Texts in Mathematics,~152. -- Berlin: Springer-Verlag, 1998.
\end{thebibliography}
\end{document}